\documentclass[12pt]{article}

\usepackage{amscd,amsmath,amssymb,eufrak}

\newtheorem{definition}{Definition}
\newtheorem{theorem}{Theorem}
\newtheorem{lemma}{Lemma}
\newtheorem{proposition}{Proposition}

\def\text#1{\mbox{#1}}

\def\Cal{\cal}

\def\frak#1{\EuFrak{#1}}

\def\la{\lambda}
\def\de{\delta}

\def\om{\omega}

\def\Si{\Sigma}
\def\al{\alpha}

\def\Th{\Theta}

\def\Ga{\Gamma}

\def\ti{\tilde}

\def\vp{\frak S^\infty}
\def\is{\frak S_\infty}
\def\Mu{{\rm M}}

\def\le{\leqslant}
\def\ge{\geqslant}

\begin{document}
\rightline{\Large{}PDMI PREPRINT --- 03/2001}
\vspace{.5in}
\begin{center}
{\Large\bf{}ON THE SIMPLEST SPLIT-MERGE OPERATOR\\[6pt]
ON THE INFINITE-DIMENSIONAL SIMPLEX} \\[0.5in]
\begin{tabular}{cc}
{\large\bf Natalia TSILEVICH\footnotemark[1]{}} \\[6pt]
St.~Petersburg State University \\
E-mail: natalia@pdmi.ras.ru & 
\end{tabular}
\end{center}
\vspace{.25in}
\centerline{\large May 31, 2001}
\vspace{.5in}
\begin{abstract}
We consider the simplest split-merge Markov operator $T$ on the
infinite-dimensional simplex $\Si_1$ of monotone non-negative sequences
with unit sum. For a sequence $x\in\Si_1$, it picks a size-biased sample (with
replacement) of two elements of $x$; if these elements are distinct, it merges
them and reorders the sequence, and if the same element is picked twice,
it splits this element
uniformly into two parts and reorders the sequence.
We prove that the means along the $T$-trajectory of the $\de$-measure at
the vector
$(1,0,0,{\ldots})$ converge to the Poisson--Dirichlet distribution $PD(1)$.
\end{abstract}

{\par\phantom{A}\par\vfill\footnoterule{}\parindent=1em\noindent
\hbox to 1.8em{\hss\normalsize${}^{1}$}\normalsize{%
Partially supported by RFBR grant 00--15--96060. 
}}

%------------------------------------------------------------------------
\newpage
\section{Introduction}

In this paper we investigate the simplest
split-merge Markov operator $T$ originally introduced by A.~Vershik.
This operator acts
on the infinite-dimensional simplex $\Si_1$
of monotone non-negative sequences with unit sum
as follows. For
$x=(x_1,x_2,\ldots)\in\Si_1$, consider a size-biased sample (with
replacement) of two elements of $x$. If two different elements $x_i$,
$x_j$ are
picked, then we merge them into $x_i+x_j$
(and reorder the sequence to obtain a point of
$\Si_1$); and if the same element $x_i$ is picked twice, then
we split it into two parts $t$, $x_i-t$
(with $t$ uniformly distributed on $[0,x_i]$) and reorder the sequence.

\smallskip\noindent{\bf Conjecture} (A.~Vershik) {\it
The only $T$-invariant distribution on $\Si_1$ is the Pois\-son--Di\-rich\-let
distribution $PD(1)$.}
\smallskip

The Poisson--Dirichlet distribution $PD(1)$ introduced by J.~F.~C.~Kingman
\cite{Ki75} is perhaps the most distinguished distribution on
the infinite-dimensional simplex. It arises
in many various problems from different areas
of mathematics and applications and
has many remarkable properties and characterizations. 
The conjecture above may be considered as an attempt to give a new
very important characterization of this measure.

In the original setting the operator $T$
arises from the representation theory of the infinite
symmetric group $\is$. Namely, the group
$\is$ acts by shifts $R_g$
on the projective limit $\vp$ of finite symmetric groups (the so-called space
of virtual permutations).
On the other hand, there is a projection from $\vp$ onto the
infinite-dimensional simplex
which associates with a virtual permutation the sequence
of its relative cycle lengths.
The simplest split-merge operator $T$ under consideration
is the projection on $\Si_1$
of the shift $R_{(1,2)}$ corresponding to a transposition in $\is$.
Following this approach,
N.~Tsilevich~\cite{Ts98, Ts99} proved the uniqueness of invariant measure for
the family of Markov operators $\{T_g\}_{g\in\is}$
arising from the action of the whole group $\is$. 
We would also like to mention the paper \cite{GK99}
proving the uniqueness of invariant measure for a closely related Markov
operator.

Recently, the same operator $T$ appeared as a simplified model
in a quite different context
of triangulations of random Riemannian surfaces (R.~Brooks). Inspired by
this motivation, E.~Mayer-Wolf, O.~Zeitouni and M.~Zerner
have proved (by purely analytic methods)
the uniqueness of $T$-invariant measure
under strong smoothness assumptions on the measure
(A.~Vershik, personal communication). 
However, the
problem for arbitrary Borel distributions is still open.

In this note we prove another partial result. Consider the trajectories
of the $\de$-measures $\de_x$, $x\in\Si_1$, under $T$. We show that 
for almost all points $x$ with finitely many non-zero coordinates (e.g.,
$x=(1,0,\ldots)$), the binomial means along
the trajectories converge to the Poisson--Dirichlet
measure $PD(1)$. 

We would like to emphasize that we
essentially exploit the relation of this problem
to symmetric groups and the techniques developped in \cite{Ts98, Ts99}. We
hope that these methods will eventually allow to prove Vershik's conjecture.

Let us describe the problem more precisely.
Let
$$
\Si=\left\{x=(x_1,x_2,{\ldots}):\;
x_1\ge x_2\ge{\ldots}\ge0,\,\sum_{i=1}^\infty x_i\le 1 \right\}
$$
be the simplex of non-increasing non-negative sequences with sum at most one, and
$\Si_1=\{x\in\Si:\,\sum_{i=1}^\infty x_i=1\}$ be its subsimplex consisting of
sequences with sum exactly one. The most elementary definition of the
Poisson--Dirichlet measure $PD(1)$ is as follows.

\begin{definition}\label{def:PD}
Let $U_1,U_2,{\ldots}$ be a sequence of i.i.d.~random variables uniformly
distributed on the interval $[0,1]$. Let
$$
V_n=U_n\,\prod_{i=1}^{n-1}(1-U_i).
$$
It is easy to see that $\sum_{i=1}^\infty V_i=1$.
The Poisson--Dirichlet measure $PD(1)$ on the simplex $\Si_1$
is the distribution of the order
statistics $V_{(1)}\ge V_{(2)}\ge{\ldots}$
of the sequence $V_1,V_2,{\ldots}$.
\end{definition}

Now consider the following Markov operator on $\Si_1$:
\begin{equation}\label{operator}
Tx\!=\!\begin{cases}
V(x_i+x_j,x_1,{\ldots},x_{i-1},x_{i+1},{\ldots},x_{j-1},x_{j+1},{\ldots}),
\quad i<j,
\\ \qquad\quad
\qquad\qquad\qquad\mbox{ with probability } 2x_ix_j;\\
V(t,x_i-t,x_1,{\ldots},x_{i-1},x_{i+1},{\ldots}),\quad i=1,2,{\ldots},
\\ \qquad\quad
\qquad\qquad\qquad\mbox{ with probability } x_idt,\,t\in[0,x_i],
\end{cases}
\end{equation}
where the operator $V$ arranges the elements of a sequence
in non-increasing order.
Denote by $E$ the identity operator on $\Si_1$. We may consider these operators
also on the space ${\cal M}(\Si_1)$
of probability Borel measures on $\Si_1$.

Let $\nu_0$ denote the $\de$-measure at the point
$(1,0,0,\ldots)\in\Si_1$. In this note we prove the following theorem.

\begin{theorem}\label{th:our}
The sequence
$$
\left(\frac{E+T}2\right)^m\nu_0
$$
converges (weakly) as $m\to\infty$ to
the Poisson--Dirichlet measure $PD(1)$.
\end{theorem}

\smallskip\noindent{\bf Corollary}
{\it The sequence
$$
\left(\frac{E+T}2\right)^m\de_x
$$
converges (weakly) as $m\to\infty$ to
the Poisson--Dirichlet measure $PD(1)$ for almost all (with respect to
the Lebesgue measure) points $x\in\Si_1$ 
with finitely many non-zero coordinates.}
\smallskip

The author is grateful to A.~Vershik for suggesting the problem and  for many 
fruitful discussions. 

\section{The space of virtual permutations}

In this Section we give a necessary background concerning the space of
virtual permutations and central measures on this space.

\begin{definition}\label{def:induced}
Given a subset $J \subset [n]$ and a permutation $w \in
{\frak S}_n$, denote by $\pi_{n,J} w$ a permutation of the set
$J$ obtained by removing from cycles of $w$ all elements that do not
belong to $J$. The permutation
$\pi_{n,J} w$ is called {\it the induced permutation} for $w$ on $J$.
\end{definition}

The induced permutation on the subset $J=[m]$ is denoted by 
$\pi_{n,m} w$. We shall usually omit the index $n$ if it is clear from the
context. For example, if
$w=(6\,3\,5\,1)\,(4\,2)\,(7)$, then 
$\pi_4 w=(3\,1)\,(4\,2)$.

\begin{definition}[\cite{KOV93}]\label{def:vp}
{\it The space of virtual permutations} ${\frak S}^\infty$
is the projective limit ${\frak S}^\infty=\varprojlim {\frak S}_n$ 
of finite symmetric groups ${\frak S}_n$ with respect to canonical
projections
$\pi_{n+1,n}:{\frak S}_{n+1}\to {\frak S}_n$.
\end{definition}

Thus a virtual permutation is a sequence
$(w_1,w_2,\,\ldots) \in 
{\frak S}_1\times {\frak S}_2\times{\ldots}$,
such that $\pi_nw_{n+1}=w_n$ for all $n\in\mathbb N$.

Note that $\pi_n$ commutes with shifts on elements of the group
${\frak S}_n$, i.e.~for all
$N>n$, the equality
$\pi_n(g_1^{-1}hg_2)=g_1^{-1}\pi_n(h)g_2$
holds for all $h\in {\frak S}_N$ and $g_1,g_2\in {\frak S}_n$.
Let $\is=\cup_{n\ge1}{\frak S}_n$ be the infinite symmetric group.
Then the group
$G={\frak S}_\infty\times {\frak S}_\infty$ acts on the space of virtual
permutations ${\frak S}^\infty$ as
\begin{equation}
((g_1,g_2)\om)_i=\begin{cases} g_1^{-1}w_ig_2, &\text{ if } i\ge n;\\
\pi_i(g_1^{-1}w_n g_2), &\text{ if } i<n.
\end{cases}
\label{action}
\end{equation}

A sequence $\{\mu_n\}$ of distributions on finite symmetric groups
${\frak S}_n$ is called {\it coherent}, if it accords with projections
$\pi_{n+1,n}$, 
i.e.~$\pi_{n+1,n}\mu_{n+1}=\mu_n$ for all $n\in\mathbb N$.
Each coherent sequence of distributions defines a Borel measure
$\mu=\varprojlim\mu_n$ on the space of virtual permutations ${\frak S}^\infty$,
and each Borel measure on $\vp$ can be presented in this form. {\it In what
follows all measures are assumed to be Borel and normalized.}
If all measures $\mu_n$ are central, i.e.~invariant under inner
automorphisms of ${\frak S}_n$, then the limit distribution
$\mu$ is invariant under the diagonal subgroup
$K=\{(g_1,g_2)\in G:\;g_1=g_2\}$. 

\begin{definition}
$K$-invariant
measures on the space of virtual permutations
${\frak S}^\infty$ are called
{\it central}. The set of all central measures is denoted by
${\Cal M}^K({\frak S}^\infty)$.
\end{definition}

\smallskip\noindent
{\bf Example 1.}
Let $m_n$ be the Haar measure on the symmetric group ${\frak S}_n$. The
sequence $\{m_n\}$ is coherent, and the measure
$m=\varprojlim m_n$ is called the {\it Haar measure 
on the space of virtual permutations}.
It is clear that
this measure is central. Moreover, it is invariant under
the whole group $G$.

The problem of describing all central measures on
the space of virtual permutations is parallel to Kingman's theory
of {\it partition structures} (see, e.g.,~\cite{Ki75, Ki78}).
Let us reformulate his results in our terms.

Each point $x\in\Si$ defines a central measure $P^x$ on the space of
virtual permutations. For $x\in\Si_1$ they are described as follows.
Let us put elements
$1,2,{\ldots}$ at random into cycles and label at random these cycles
according to the following rule:

\begin{itemize}
\item at the first step the element~$1$ forms a $1$-cycle, and we 
label it by $j$ with probability
$x_j$ ($j=1,2,{\ldots}$);
\item if the elements $1,{\ldots},m$ are already placed, and they form
$k$ cycles with labels $i_1,{\ldots},i_k$,
then the element
$m+1$ is inserted in one of possible positions in the $j$th
cycle with probability $x_{i_j}$ ($j=1,{\ldots},k$), or forms a new
$1$-cycle labelled by $i$ with probability
$x_i$ ($i\ne i_1,{\ldots},i_k$).  
\end{itemize}
After the $n$th step of this procedure we obtain a random permutation
$w_n\in {\frak S}_n$. It follows from construction that the sequence
$\{w_n\}$ is coherent, thus it defines a random element of $\frak S^\infty$,
and we denote its distribution by $P^x$.

\smallskip\noindent{\bf Example 2.}
If $x_0=(1,0,{\ldots})$, then the projection of the measure
$\tau=P^{x_0}$ onto $\frak S_n$
is the measure $\tau_n$ uniformly distributed on
one-cycle permutations.

\smallskip\noindent
{\bf Example 3.} The projection of the measure
$P^x$, $x\in\Si_1$, on ${\frak S}_2$ is given by
\begin{eqnarray*}
P^x_2((1)(2))&=&2\sum_{1\le i<j<\infty} x_ix_j,\\
P^x_2((12))&=&\sum_{i=1}^\infty x_i^2.
\end{eqnarray*}

\begin{theorem}[\cite{Ki78}]\label{th:Kingman}
Let $\mu$ be a central measure on the space of virtual permutations.
Denote by $l_1(w_n)\ge l_2(w_n)\ge{\ldots}$
the cycle lenghts of a permutation $w_n\in \frak S_n$ in non-increasing order.
Then the limits
\begin{equation}
\label{2}
X_i(\om)=\lim_{n\to\infty}\frac{l_i(w_n)}n,\qquad i=1,2,{\ldots},
\end{equation}
called the relative cycle lengths of the virtual permutation $\om$,
exist for almost all with respect to $\mu$ virtual permutations
$\om=(w_1,w_2,\ldots)$.
For every $x\in\Si$,
the conditional distribution of $\om$, given
$(X_1,X_2,{\ldots})=x$, equals $P^x$. Thus
\begin{equation}
\mu=\int_\Si P^x d\nu(x),
\label{3}
\end{equation}
where $\nu$ is the distribution of the vector $X=(X_1,X_2,{\ldots})$
on $\Si$.
\end{theorem}

In particular, if the measure $\mu$
is {\it saturated}, i.e.~the sum of relative cycle lengths
is equal to 1 a.s.~with respect to $\mu$, then the corresponding
measure $\nu$ is supported by $\Si_1$. 

\smallskip\noindent
{\bf Example 4.}
The distribution of the relative cycle lengths of virtual
permutations with respect to the 
Haar measure is the
Poisson--Dirichlet distribution $PD(1)$.
\smallskip

Thus there is a one-to-one correspondence $\rho:\nu\mapsto\mu$
between Borel distributions $\nu$
on $\Si$ and central distributions $\mu$ on $\vp$.
Given $\nu$, the
corresponding central measure $\mu=\rho(\mu)$ is recovered via~(\ref{3}).
And given $\mu$, the
corresponding measure $\nu=\rho^{-1}(\mu)$ on $\Si$
is the distribution of the relative cycle
lengths with respect to $\mu$. 

\smallskip\noindent
{\bf Remark.} The space of virtual permutations $\vp$ was introduced 
by S.~Kerov, G.~Olshanski and A.~Vershik
\cite{KOV93}. They considered a family of quasi-invariant 
distributions on $\vp$ with respect to the action of the group
$G=\is\times\is$ and studied the associated family of unitary representations
of the infinite symmetric group.

\section{Reduction lemmas}

The infinite symmetric group ${\frak S}_\infty$ acts on
the space of virtual permutations
by two-sided shifts~(\ref{action}).
Denote by $R_g:\vp\to\vp$ the right shift on
$g\in {\frak S}_\infty$.

One can easily check that the image $\mu^g=R_g\mu$ of a
central
measure $\mu$ enjoys the property that the normalized cycle lengths
exist for almost all with respect to $\mu^g$
virtual permutations. Hence we may define the projection of $R_g$ on
$\Si$, which is already a Markov operator.
It sends a point $x\in\Si$ to $\rho^{-1}(\om g)$, where
$\om$ is a random virtual permutation with distribution $P^x$.
Thus the corresponding operator $T_g$ on the space
$M(\Si)$ of probability measures on $\Si$
completes the commutative diagram
$$
\begin{CD}
M^K({\frak S}^\infty) @> R_g >> M({\frak S}^\infty) \\
@A \rho AA  @V \rho^{-1}  VV \\
M(\Si) @>T_g >> M(\Si).
\end{CD}
$$
One can easily see that the operator $T_g$ depends in fact only on
the conjugacy class of $g$, i.e.~on its cycle structure.

The subsimplex
$\Si_a=\{x\in\Si:\,\sum x_i=a\}$
is $T_g$-invariant
for all $a\in[0,1]$ and all $g\in\is$ (since
the shift by a finite permutation does not change the sum of the
relative cycle lengths of a virtual permutation).
It is proved in \cite{Ts99} that
these sets are ergodic components, and the ergodic measure concentrated
on $\Si_a$ is the image of $PD(1)$ under the homothetic transformation
$\Ga_a:\Si_1\to\Si_a$. In particular,

\begin{theorem}[\cite{Ts99}] 
The only measure on $\Si_1$ which is invariant
under the family $\{T_g\}_{g\in {\frak S}_\infty}$ is the Poisson--Dirichlet
distribution $PD(1)$.
\end{theorem}

The relation of this theorem to our problem
is given by the following lemma (note that $T_e$ is obviously the identity
operator, and we may consider
the identity operator $E$ on $\Si_1$ as
the restriction of $T_e$ on $\Si_1$). 

\begin{lemma}
The Markov operator $T$ on $\Si_1$ defined by~(\ref{operator})
is the restriction on $\Si_1$
of the operator $T_{(1,2)}$ ($=T_g$, where $g$ is an arbitrary
transposition in $\is$).
\end{lemma}

\smallskip\noindent{\it Proof.} Follows from the definition of $T_g$
and Example~3.

Now we cite some preliminary
lemmas from \cite{Ts99} which we shall use in the sequel.
First, we reformulate the problem in terms of 
the space of virtual permutations. Since there is a one-to-one
correspondence between Borel measures on $\Si$ and central measures
on $\vp$, one may think of 
the operator $T_g$ as acting on the space
${\cal M}^K(\vp)$. More exactly, denote by $\ti T_g$ the operator on
${\cal M}^K(\vp)$, corresponding to $T_g$, i.e.,
the image of $T_g$ under the map $\rho$.  

\begin{lemma}[\cite{Ts99}]\label{l:P}
Let $\mu\in{\cal M}(\vp)$. Then
$\ti T_g\mu=PR_g\mu$, where $P:M(\frak S^\infty)\to
M^K(\frak S^\infty)$ is the projection (conditional expectation)
onto the space of central measures. 
\end{lemma}

In particular, let $\ti T=\ti T_{(1,2)}$ and
$\ti E=\ti T_e$ (the identity operator on ${\cal M}(\vp)$).
Then our Theorem~\ref{th:our} is equivalent to the following proposition.
Let $\tau=\varprojlim\tau_n=P^{x_0}$ 
be the central measure on $\vp$ corresponding to the  
vector $x_0=(1,0,0,{\ldots})\in\Si_1$ (see Example~3).

\begin{proposition}\label{pr:our}
Let $\tau^q=\left(\frac{\ti E+\ti T}2\right)^q\tau$. Then
$$
\lim_{q\to\infty}\tau^q=m,
$$
where $m$ is the Haar measure on the space of virtual permutations.
\end{proposition}

The ergodic method allows one to express the operator $\ti T_g$
in the following form.
\begin{lemma}[\cite{Ts99}]\label{l:erg}
Let $\mu$ be a central measure on the space of virtual permutations. Then
\begin{equation}
\ti T_g\mu=\lim_{N\to\infty}\frac1{N!}\sum_{u\in \frak S_N}u^{-1}(R_g\mu)u.
\label{eq:proj}
\end{equation}
\end{lemma}

Finally, the formula~(\ref{eq:proj}) may be restated in terms of
distributions on {\it finite} symmetric groups.

\begin{lemma}[\cite{Ts99}]\label{l:main}
Let $\mu$ be a central measure on the space of virtual permutations, and
$\mu_n$ be the corresponding coherent family of measures
on finite symmetric groups. Then the finite-dimensional projections
of the measure $\ti T_g\mu$ are given by the following formula,
\begin{equation}
(\ti T_g\mu)_n(u)=
\sum_{w\in {\frak S}_{n+k}\atop\pi_n(w)=u}\mu_{n+k}(wg^n)\qquad
\forall n,k\in\Bbb N,\;\forall g\in\frak S^k,\;\forall u\in {\frak S}_n,
\label{eq:cond}
\end{equation}
where ${\frak S}^k={\frak S}[n+1,{\ldots},n+k]\subset {\frak S}_{n+k}$.
\end{lemma}

\section{Shifted projection of irreducible characters
of symmetric groups}
\label{sec:char}

In this section we present a formula obtained in \cite{Ts98}
which describes the action of the operator $\ti T_g$
on irreducible characters $\chi_\la$ of the symmetric group
$\frak S_{n+1}$. It generalizes a known formula for the action
on $\chi_\la$ of the canonical projection $\pi_n$.

Fix $n,k\in\Bbb N$, $N=n+k$ and
$g\in\frak S^k=\frak S[n+1,{\ldots},n+k]$.
Define a projection $\pi^g_{N,n}:\frak S_N\to \frak S_n$ as 
$$
\pi^g_{N,n}(h)=\pi_n(hg).
$$
Then the formula~(\ref{eq:cond}) from
Lemma~\ref{l:main} may be formulated as follows: 
$$
(\ti T_g\mu)_n=\pi^g_{N,n}\mu_N
$$
for all $n,k\in\Bbb N$, $N=n+k$, $g\in\frak S^k$.
We may consider the central measure $\mu_n$ as a central function on
the symmetric group $\frak S_n$. Since the characters of irreducible representations
form a basis in the space of central functions, it is interesting to consider
the action of $\pi^g_{N,n}$ on the irreducible characters. Recall that the
irreducible representations of the symmetric group $\frak S_n$ are indexed by
partitions $\la=(\la_1\ge\la_2\ge\ldots)$
of the number $n$ (which we identify with the corresponding Young
diagrams). 
Denote by $\Pi_n$ the set of all
partitions of $n$. Thus we are interested in a formula for
$$
\pi^g_{N,n}\chi_\la(u)=\sum_{w\in \frak S_N\atop \pi_n(w)=u}\chi_\la(wg),
\quad u\in \frak S_n.
$$

\smallskip\noindent {\bf Remark.} 
If $g=e$ this is the action of the canonical projection 
$\pi_n$ on the irreducible characters. As follows from the arguments of
\cite[Sect.~4.6]{Ol89},
\begin{equation}
\pi_n\chi_\la=
\sum_{\mu\nearrow\la}
\left(c(\la\setminus\mu)+1\right)\chi_\mu,
\qquad \la\in\Pi_{n+1},
\label{eq:chi}
\end{equation}
where $\mu\nearrow\la$ means that the diagram 
$\mu\in\Pi_n$ is obtained from
$\la$ by removing one cell. 

We recall some definitions related to Young diagrams.
The {\it Young diagram } $D_\la$ of the partition $\la\in\Pi_n$
is the set
$\{(i,j):\,i,j\in\Bbb Z,\, i\ge1,\, \la_i\ge j\ge1\}$. If $x=(i,j)\in D_\la$,
then $x$ is called a {\it cell} of the diagram $D_\la$.
The {\it content} of a cell $x=(i,j)$
is the number $c(x)=j-i$.
Let $m<n$, $\la\in\Pi_n$,
$\mu\in\Pi_m$. We write $\la\supset\mu$ if
$\la_i\ge\mu_i$ for all $i=1,2,{\ldots}$, that is, the Young diagram of
$\la$ contains the Young diagram of
$\mu$. The set-theoretic difference of these diagrams is called a
{\it skew diagram}. A skew diagram is called a
{\it skew hook}, if it is connected and does not contain two cells
in one diagonal. If a skew hook contains
$k$ cells, it is called a
{\it skew $k$-hook}. Denote by $\Th_k$ the set of skew
$k$-hooks.
The {\it height} of a skew hook $\la\setminus\mu$
is the number $l(\la\setminus\mu)$ equal to the number of rows it intersects minus
one. A skew diagram is called a {\it horizontal $m$-strip},
if it contains $m$ cells and
has at most one cell in each row. 

\begin{proposition}[\cite{Ts98}]\label{prop:skew}
Let $\chi_\la$ be the character of the irreducible representation
of the symmetric group $\frak S_N$ corresponding to the partition
$\la\in\Pi_N$, and let $g\in\frak S^k$ be a permutation consisting of one 
cycle of length $k$. Then
\begin{equation}
\pi^g_{N,n}\chi_\la=\sum_
{ \la\setminus\mu\in\Th_k}
(-1)^{l(\la\setminus\mu)}
\prod_{x\in\la\setminus\mu}\left(c(x)+1\right)\chi_\mu.
\end{equation}
\end{proposition}

\medskip\noindent
{\bf Corollary} {\it 
Let $\chi_\la$ be the character of the irreducible representation of
the symmetric group $\frak S_n$ corresponding to the partition
$\la\in \Pi_{n+k}$, and let $g\in \frak S^k$ be a permutation with cycle  
structure $\nu=(\nu_1,{\ldots},\nu_s)\in\Pi_k$. Then
\begin{equation}
\pi^g_{N,n}\chi_\la=\sum_\Mu
(-1)^{l(\Mu)}
\prod_{x\in\la\setminus\mu}(c(x)+1)\chi_\mu,
\end{equation}
where the sum is taken over all $s$-tuples $\Mu=\{\mu_1,{\ldots},\mu_s\}$
such that 
$\la=\mu_0\supset\mu_1\supset{\ldots}\supset\mu_s\in\Pi_n$, 
$\mu_{i-1}\setminus\mu_{i}$ is a skew $\nu_i$-hook for all
$i=1,{\ldots},s$, and
$l(\Mu)=\sum_{i=1}^s
l(\mu_{i-1}\setminus\mu_{i})$.}

\section{Proof of Proposition~\ref{pr:our}}

Let $N=n+2m$, and $\la\in\Pi_N$.

\begin{lemma}
\label{l:calc}
\begin{equation}
\frac1{2^m}(\ti E+\ti T)^m\chi_\la=
\sum_{\mu\in\Pi_n}\chi_\mu\cdot 
\left(\prod_{x\in\la\setminus\mu}(c(x)+1)\right)\cdot
p(\la,\mu),
\end{equation}
where $p(\la,\mu)$ is the number of all sequences
$\Mu=\{\mu_0=\mu,\mu_1,{\ldots},\mu_m=\la\}$ such that
$\mu_k\in\Pi_{n+2k}$, and $\mu_{k+1}\setminus\mu_k$ is
a horizontal $2$-strip for all $k=0,\ldots,m-1$.
\end{lemma}

\smallskip\noindent{\it Proof.}
Note that
\begin{multline*}
\frac1{2^m}(\ti E+\ti T)^m\chi_\la=\\
\left(\frac{\pi_{n+2,n}+\pi_{n+2,n}^{\rm{tr}}}{2}\right)
\left(\frac{\pi_{n+4,n+2}+\pi_{n+4,n+2}^{\rm{tr}}}{2}\right)
{\ldots}
\left(\frac{\pi_{N,N-2}+\pi_{N,N-2}^{\rm{tr}}}{2}\right)
\chi_\la,
\end{multline*}
where $\pi_{l+2,l}^{\rm{tr}}= \pi_{l+2,l}^{(l+1,l+2)}$. 
Then
it follows from Proposition~\ref{prop:skew} that the projection of
$\frac1{2^m}(E+T)^m\chi_\la$ onto ${\frak S}_n$
equals
\begin{equation*}
\sum_{\mu\in\Pi_n}\chi_\mu\cdot 
\left(\prod_{x\in\la\setminus\mu}(c(x)+1)\right)\cdot 
\sum_\Mu \al(\Mu),
\end{equation*}
where the sum is taken over all sequences
$\Mu=\{\mu_0=\mu,\mu_1,{\ldots},\mu_m=\la\}$ such that
$\mu_{k}\subset\mu_{k+1}$, $\mu_k\in\Pi_{n+2k}$, and
$\al(\Mu)=\prod_{k=1}^m\al(\mu_k\setminus\mu_{k-1})$, where
$\al(\cdot)=\frac{\al_1(\cdot)+\al_2(\cdot)}{2}$ with
\begin{eqnarray*}
\al_1(\la\setminus\nu)&=&\cases
1, &\text{if } \la\setminus\nu \text{ is a 2-row or a 2-column},\\
2, &\text{otherwise}.
\endcases\\
\al_2(\la\setminus\nu)&=&\cases
1, &\text{if } \la\setminus\nu \text{ is a 2-row},\\
-1, &\text{if } \la\setminus\nu \text{ is a 2-column},\\
0, &\text{otherwise},
\endcases
\end{eqnarray*}
that is,
$$
\al(\la\setminus\nu)=\cases
1, &\text{if } \la\setminus\nu \text{ is a horizontal 2-strip},\\
0, &\text{otherwise (that is, if $\la\setminus\nu$ is a 2-column)}.
\endcases
$$
The Lemma follows.

Now each central measure on ${\frak S}_N$ can be decomposed into a combination of
irreducible characters. In particular, let
$$
\tau_N(g)=\sum_{\la\in\Pi_N}a_\la\chi_\la(g).
$$
It is easy to see that the coefficients of this decomposition equal
$$
a_\la=\frac1{N!}\chi_\la((12{\ldots}N)).
$$
It follows from the Murnaghan--Nakayama rule that these coefficients are not zero
only for hook partitions $\la_k^N=(k1^{N-k})$, for $k=1,{\ldots},N$, and
$a_{(k1^{N-k})}=\frac{(-1)^{N-k}}{N!}$. Thus
$$
\tau_N=\frac1{N!}\sum_{k=1}^N (-1)^{N-k}\chi_{(k1^{N-k})}.
$$
Now we are ready to prove Proposition~\ref{pr:our}.

First let $n=2$. Denote $c(\la,\mu)=\prod_{x\in\la\setminus\mu}(c(x)+1)$.
Then the projection $\tau_2^q$ of $\tau^q$ on ${\frak S}_2$ can be
written as $a_q\chi_{(2)}+b_q\chi_{(1^2)}$, where
\begin{eqnarray*}
a_q&=&\sum_{k=1}^N\frac{(-1)^k}{N!}\cdot c(\la_k^N,(2))\cdot p(\la_k^N,(2)),\\
b_q&=&\sum_{k=1}^N\frac{(-1)^k}{N!}\cdot c(\la_k^N,(1^2))\cdot 
p(\la_k^N,(1^2)).
\end{eqnarray*}
It is clear that the coefficient
$p(\la_k^N,(2))$ is not zero
only for $k=N$, and
$c((N),(2))=3\cdot4\cdot{\ldots}\cdot N=N!/2$. Thus we obtain that
$a_q=1/2$. The corresponding
term $\frac12\chi_{(2)}=m_2$ is just the Haar measure on ${\frak S}_2$.
Thus it suffices to show that $b_q\to0$ as $q\to\infty$.

Let us calculate the coefficient $p(\la_k^N,(1^2))$. It follows from
Lemma~\ref{l:calc}
 that we must calculate the number of paths from $(1^2)$ to $(k1^{N-k})$ such
that at each step we add either two cells in one row or two cells in different
rows and different columns.
It is easy to see that this number
equals $q\choose{N-k-1}$. The content coefficient equals
$c((k1^{N-k}),(1^2))=(-1)^{k+1}k!(N-k-1)!$, thus the summand in $b_q$
corresponding to $\la_k^N$ equals
$$
-\frac{q!k!(N-k-1)!}{N!(N-k-1)!(k-q-1)!}=
-\frac{k!q!}{N!(k-q-1)!}.
$$
Hence,
\begin{multline*}
b_q=\frac{q!}{N!}\sum_{k=q+1}^{2q+1}\frac{k!}{(k-q-1)!}=
\frac{q!(q+1)!}{N!}\sum_{k=q+1}^{2q+1}\frac{k!}{(q+1)!(k-q-1)!}=\\
\frac{q!(q+1)!}{N!}\sum_{k=q+1}^{2q+1}{k\choose{q+1}}=
\frac{q!(q+1)!}{N!}{{2q+2}\choose{q+2}}=\frac{1}{q+2}\to0
\end{multline*}
as $q\to\infty$, as desired.

The proof for $n>2$ goes in the same way with a bit complicated calculations.
It is clear that the decomposition of $\tau^q_n$ contains only hook characters,
that is,
$$
\tau^q_n=\sum_{l=1}^n a_q(l)\chi_{(l1^{n-l})},
$$
and
$$
a_q(l)=
\sum_{k=1}^N\frac{(-1)^{N-k}}{N!}\cdot c(\la_k^N,\la_l^n)\cdot p(\la_k^N,\la_l^n).
$$

If $l=n$, there is again only one non-zero summand corresponding to $k=N$, 
$$
a_q(n)=\frac{1}{N!}\cdot\frac{N!}{n!}\cdot 1=\frac1{n!},
$$
and $a_q(n)\chi_{(n)}=m_n$, that is, the Haar measure on ${\frak S}_n$. Thus
we only have to prove that $a_q(l)\to0$ as $q\to\infty$ for $l<n$.
Easy calculations show that
$p(\la_k^N,\la_l^n)={q\choose{2q-k+l}}$, and
$c(\la_k^N,\la_l^n)=(-1)^{l-k}\frac{k!(N-k-1)!}{l!(n-l-1)!}$.
Thus 
$$
a_q(l)=\frac{(-1)^{l+N}}{l!{(n-l-1)!}}\cdot\frac{q!}{N!}\cdot\sum_{k=q+l}^{2q+l}
\frac{k!(N-k-1)!}{(2q-k+l)!(k-l-q)!}.
$$
But
$\frac{(N-k-1)!}{(N-k-n+l)!}\le(N-k-1)^{n-l-1}\le(q+n-l-1)^{n-l-1}$. Then
\begin{eqnarray*}
|a_q(l)|&\le&{\rm const}\cdot\frac{(q+n-l-1)^{n-l-1}q!(q+l)!}{N!}
\sum_{k=q+l}^{2q+l}{k\choose{q+l}}\le\\
&\le&{\rm const}\cdot\frac{(q+n-l-1)^{n-l-1}q!(q+l)!}{N!}\cdot
{{2q+l+1}\choose{q+l+1}}\le\\
&\le&{\rm const}\cdot\frac{(q+n-l-1)^{n-l-1}}{q+l+1}\cdot\frac{(2q+l+1)!}{(2q+n)!}\le\\
&\le&{\rm const}\cdot\frac{1}{q+l+1}\left(\frac{q+n-l-1}{2q+l+2}\right)^{n-l-1}\to0
\end{eqnarray*}
as $q\to\infty$. Proposition follows.

\end{document}